\documentclass[11pt,letterpaper]{article}
\usepackage{graphicx}
\usepackage{amsthm}
\usepackage{amssymb}

\usepackage[margin=1in]{geometry}

\usepackage{mathtools}

\newtheorem{theorem}{Theorem}[section]
\newtheorem{lemma}[theorem]{Lemma}

\newtheorem{remark}[theorem]{Remark}
\newtheorem{example}[theorem]{Example}

\graphicspath{{figures/}}

\newcommand{\R}{{\mathbb R}}

\title{Approximation of the joint spectral radius \\ using sum of squares}
\author{Pablo A. Parrilo\thanks{Laboratory for Information and Decision Systems, Massachusetts Institute of Technology, \texttt{parrilo@mit.edu}}
\and
Ali Jadbabaie\thanks{GRASP Laboratory, University of Pennsylvania, \texttt{jadbabai@seas.upenn.edu}}}

\date{\today}
\date{}

\begin{document}
\maketitle

\begin{abstract}
We provide an asymptotically tight, computationally efficient
approximation of the joint spectral radius of a set of matrices using
sum of squares (SOS) programming. The approach is based on a search
for an SOS polynomial that proves simultaneous contractibility of a
finite set of matrices. We provide a bound on the quality of the
approximation that unifies several earlier results and is independent
of the number of matrices. Additionally, we present a comparison
between our approximation scheme and earlier techniques, including the
use of common quadratic Lyapunov functions and a method based on
matrix liftings. Theoretical results and numerical investigations show
that our approach yields tighter approximations.
\end{abstract}

\section{Introduction}

Stability of discrete linear inclusions has been a topic of major
research over the past two decades. Such systems can be represented as
a switched linear system of the form $x(k+1) = A_{\sigma(k)} x(k)$,
where $\sigma$ is a mapping from the integers to a given set of
indices. The above model, and its many variations, has been studied
extensively across multiple disciplines including control theory,
theory of non-negative matrices and Markov chains, subdivision schemes
and wavelet theory, dynamical systems, etc. The fundamental question
of interest is to determine whether $x(k)$ converges to a limit, or
equivalently, whether the infinite matrix products chosen from the set
of matrices converge~\cite{BeWa92,DaLa92,DaLa01}.  The research on
convergence of infinite products of matrices spans across four
decades. A majority of results in this area has been provided in the
special case of non-negative and/or stochastic matrices. A
non-exhaustive list of related research providing several necessary
and sufficient conditions for convergence of infinite products and
their applications
includes~\cite{CH94,DaLa01,Leiz92,ShuWuPa97}. Despite the wealth of
research in this area, finding algorithms that can unambiguously
decide convergence remains elusive.  Much of the difficulty of this
problem stems from the hardness in computation or efficient
approximation of the joint spectral radius of a finite set of
matrices.  This notion was introduced by Rota and Strang \cite{RoSt60}
via the definition
\begin{equation}
\rho(A_1,\ldots,A_m) := \lim_{k \rightarrow \infty}
\max_{\sigma \in \{1,\ldots,m\}^k} || A_{\sigma_k} \cdots
A_{\sigma_2} A_{\sigma_1} ||^{1/k},
\label{eq:defjsr}
\end{equation}
and represents the maximum growth rate that can be achieved by taking
arbitrary products of the matrices $A_i$. As in the case of the
classical spectral radius, the value of this expression is independent
of the choice of norm in~(\ref{eq:defjsr}). Daubechies and
Lagarias~\cite{DaLa92} conjectured that the joint spectral radius is
equal to a related quantity, the {\it generalized spectral radius},
which is defined in a similar way except for the fact that the norm of
the product is replaced by the spectral radius. Berger and
Wang~\cite{BeWa92} proved this conjecture to be true for finite sets
of matrices.  Blondel and Tsitsiklis have shown that computing $\rho$
is hard from a computational complexity viewpoint, and even
approximating it is difficult~\cite{BlTi2,BlTi3}. In particular, it
follows from their results that the problem ``Is $\rho \leq 1$?'' is
undecidable. For rational matrices, the joint spectral radius is not a
semialgebraic function of the data, thus ruling out a very large class
of methods for its exact computation. We refer the reader to the
survey \cite[\S3.5]{BlTi1} for further results and references on the
computational complexity of the joint spectral radius.

It turns out that a necessary and sufficient condition for the
stability of a linear difference inclusion is for the corresponding
matrices to have a subunit joint spectral radius, i.e.,
$\rho(A_1,\ldots,A_m) < 1$; see e.g. \cite[Thm.~1]{ShuWuPa97} and
\cite{BraytonTong2}. A subunit joint spectral radius is equivalent to
the existence of a common norm with respect to which all matrices in
the set are contractive~\cite{Bar88,Koz90,wirth}; unfortunately, this
common norm is in general not finitely constructible. In fact a
similar result, due to Dayawansa and Martin~\cite{DayaMar}, holds for
nonlinear systems that undergo switching. A popular approach towards
approximating the joint spectral radius or showing that it is indeed
subunit has been to try to prove simultaneous contractibility (i.e.,
existence of a common norm with respect to which matrices are
contractive), by searching for a common ellipsoidal norm, or
equivalently, a common quadratic Lyapunov function. The benefit of
this approach is due to the fact that the search for a common
ellipsoidal norm can be posed as a semidefinite program and solved
efficiently using interior point techniques. However, it is not too
difficult to generate examples where the discrete inclusion is {\it
absolutely asymptotically stable}, i.e., asymptotically stable for all
switching sequences, but a common quadratic Lyapunov function, (or
equivalently a common ellipsoidal norm) does not exist.

Ando and Shih describe in~\cite{Ando98} a constructive procedure for
generating a set of $m$ matrices whose joint spectral radius is equal
to $\frac{1}{\sqrt{m}}$, but for which no quadratic Lyapunov function
exists.  They prove that the interval $[0,\, \frac{1}{\sqrt{m}})$ is
effectively the ``optimal" range for the joint spectral radius
necessary to guarantee simultaneous contractibility under an
ellipsoidal norm for a finite collection of $m$ matrices. The range is
denoted as optimal since it is the largest subset of $[0,1)$ for which
if the joint spectral radius is in this subset the collection of
matrices is simultaneously contractible under an ellipsoidal
norm. Furthermore, they show that the optimal joint spectral radius
range for a {\it bounded} set of $n \times n$ matrices is the interval
$[0,\, \frac{1}{\sqrt{n}})$. The proof of this fact is based on John's
ellipsoid theorem \cite{JohnEllipsoid}.  Roughly speaking, John's
ellipsoid theorem implies that every convex body in $n$-dimensional
Euclidean space that is symmetric with respect to the origin can be
approximated by inner and outer ellipsoids, up to a factor of
$\frac{1}{\sqrt{n}}$. Independently, Blondel, Nesterov and Theys
\cite{BlNT04} showed a similar result (also based on John's ellipsoid
theorem), that the best ellipsoidal norm approximation of the joint
spectral radius provides a lower bound and an upper bound on the
actual value. Given a set ${\mathcal M}$ of $n \times n$ matrices with
joint spectral radius $\rho$, and best ellipsoidal norm approximation
$\hat \rho$, it is shown there that
\begin{equation}
\frac{1}{\sqrt{n}} \, \hat \rho({\mathcal M}) \le \rho({\mathcal M})
\le \hat \rho({\mathcal M}).
\label{eq:sqrtn}
\end{equation}
A major consequence of these results is that finding a common Lyapunov
function becomes increasingly hard as the dimension goes up.

There have been a number of earlier works proposing different
numerical techniques for the effective computation of bounds on the
joint spectral radius.  A natural class of lower bounds is obtained by
considering periodic switching sequences, in which case only a finite
number of matrix norms need to be computed.  Using a naive approach,
the required computational efforts grow exponentially as $m^k$, where
$k$ is the period of the sequence.  Due to the cyclic property of the
spectral radius, some terms are redundant, and Maesumi \cite{Maesumi}
has shown using combinatorial techniques that the number of required
products can be reduced to $m^k/k$. Another approach is the work of
Gripenberg \cite{Gripenberg}, who has introduced a branch-and-bound
algorithm to produce upper and lower bounds on the joint spectral
radius.  Protasov \cite{Protasov1,Protasov2} has developed a geometric
method to approximate this quantity, based on a polytopic
approximation of a convex set that is invariant under the action of
the linear operators $A_i$.  This method has also been extended to the
computation of the so-called $p$-radius \cite{Protasov1}. More
recently, Blondel and Nesterov \cite{BlNes05} have proposed an
alternative scheme to the computation of the joint spectral radius, by
``lifting'' the matrices using Kronecker products to provide better
approximations.  A common feature in many of these approaches is the
presence of convexity-based methods to provide certificates of the
desired system properties.

In this paper, we develop a sum of squares (SOS) based scheme for the
approximation of the joint spectral radius. The method computes, using
the techniques of semidefinite programming, a homogeneous polynomial
that serves as a Lyapunov-like function for the corresponding switched
linear system. We prove several results on the quality of
approximation of the proposed scheme. In particular, it will follow
from Theorems~\ref{thm:sos2dbound} and~\ref{thm:msos2dbound} that our
SOS-based approximation $\rho_{SOS,2d}$ satisfies
\[
\eta^{-\frac{1}{2d}} \, \cdot \,
\rho_{SOS,2d}({\mathcal M}) \le \rho({\mathcal M})
\le \rho_{SOS,2d}({\mathcal M}),
\]
where $\eta := \min \{ m , {\textstyle\binom{n+d-1}{d}} \}$.  To prove
this, we use two different techniques, one inspired by recent results
of Barvinok~\cite{Barvinok} on approximation of norms by polynomials,
and the other one based on a convergent iteration similar to that used
for Lyapunov inequalities. Our results provide a simple and unified
derivation of most of the available bounds, including some new
ones. We prove that the SOS-based approximation is always tighter than
that obtained by the use of common quadratic Lyapunov functions, and
than the one provided by Blondel and Nesterov in
\cite{BlNes05}. Furthermore, we show how to compute the bound in
\cite{BlNes05} using matrices that are exponentially smaller than
those proposed there; this result also follows from the earlier work
of Protasov \cite{Protasov1}. A preliminary version of some of our
results has been presented in \cite{PabloAliHSCC}.

A description of the paper follows. In Section~\ref{sec:sosnorms} we
present a class of bounds on the joint spectral radius based on
simultaneous contractivity with respect to a norm, followed by a sum
of squares-based relaxation, and the corresponding suboptimality
properties. In Section~\ref{sec:symmalgebra} we present some
background material in multilinear algebra, necessary for our
developments, and a derivation of a bound of the quality of the SOS
relaxation. An alternative development is presented in
Section~\ref{sec:soslyap}, where a different bound on the performance
of the SOS relaxation is given in terms of a very natural Lyapunov
iteration, similar to the classical case. In
Section~\ref{sec:comparison} we make a comparison with earlier
techniques and analyze a numerical example. Finally, in
Section~\ref{sec:conclusions} we present our conclusions.

\section{Bounds via polynomials and sums of squares}
\label{sec:sosnorms}

A natural way of bounding the joint spectral radius is to find a
common norm that guarantees certain contractiveness properties for all
the matrices. In this section, we first revisit this characterization,
and introduce our method of using SOS relaxations to approximate this
common norm.

\paragraph{Norms and the joint spectral radius.}
As we mentioned, there exists an intimate relationship between the
spectral radius and the existence of a vector norm under which all the
matrices are simultaneously contractive. This is summarized in the
following theorem, a special case of Proposition 1 in \cite{RoSt60} by
Rota and Strang.

\begin{theorem}[\cite{RoSt60}]
\label{thm:RotaStrang}
Consider a finite set of matrices $\mathcal{A} =
\{A_1,\ldots,A_m\}$. For any $\epsilon > 0$, there exists a norm
 $\|\cdot\|$ in $\R^n$ (denoted as JSR norm hereafter) such that
\[
||A_i x|| \leq (\rho(\mathcal{A}) + \epsilon) \, ||x||, \qquad \forall x
  \in \R^n, \quad i = 1,\ldots,m.
\]
\end{theorem}

The theorem appears in this form, for instance, in Proposition~4 of
\cite{BlNT04}.  The main idea in our approach is to replace the JSR
norm that approximates the joint spectral radius with a homogeneous
SOS polynomial $p(x)$ of degree $2d$. As we will see in the next
sections, we can produce arbitrarily tight SOS approximations, while
still being able to prove a bound on the resulting estimate.

\paragraph{Joint spectral radius and polynomials.}
As the results presented above indicate, the joint spectral radius can
be characterized by finding a common norm under which all the maps are
simultaneously contractive.  As opposed to the unit ball of a norm,
the level sets of a homogeneous polynomial are not necessarily convex
(see for instance Figure~\ref{fig:jsr}). Nevertheless, as the
following theorem suggests, we can still obtain upper bounds on the
joint spectral radius by replacing norms with homogeneous polynomials.

\begin{theorem}
\label{thm:psdbound}
Let $p(x)$ be a strictly positive homogeneous polynomial of degree
$2d$ that satisfies
\[
p(A_i x) \leq \gamma^{2d} \, p(x), \qquad \forall x \in \R^n \quad i = 1,\ldots,m.
\]
Then, $\rho(A_1,\ldots,A_m) \leq \gamma$.
\end{theorem}
\begin{proof}
If $p(x)$ is strictly positive, then by compactness of the unit ball
in $\R^n$ and continuity of $p(x)$, there exist constants $0 < \alpha
\leq \beta$, such that
\[
\alpha \, ||x||^{2d} \leq p(x) \leq \beta \, ||x||^{2d} \qquad \forall x \in \R^n.
\]
Then,
\begin{eqnarray*}
||A_{\sigma_k} \ldots A_{\sigma_1}|| &\leq&
\max_x \frac{||A_{\sigma_k} \ldots A_{\sigma_1} x||}{||x||}  \\
&\leq& \left(\frac{\beta}{\alpha}\right)^\frac{1}{2d} \max_x \frac{p(A_{\sigma_k} \ldots A_{\sigma_1} x)^\frac{1}{2d}}{p(x)^\frac{1}{2d}} \\
& \leq & \left(\frac{\beta}{\alpha}\right)^\frac{1}{2d} \gamma^k.
\end{eqnarray*}
From the definition of the joint spectral radius in
equation~(\ref{eq:defjsr}), by taking $k$th roots and the limit $k
\rightarrow \infty$ we immediately have the upper bound
$\rho(A_1,\ldots,A_m) \leq \gamma$.
\end{proof}

The condition in Theorem~\ref{thm:psdbound} involves positive
polynomials, which are computationally hard to characterize.  A useful
scheme, introduced in \cite{Phd:Parrilo,sdprelax} and relatively
well-known by now, relaxes the nonnegativity constraints to a much
more tractable \emph{sum of squares} (SOS) condition, where $p(x)$ is
required to have a decomposition as $p(x) = \sum_i p_i(x)^2$. The SOS
condition can be equivalently expressed in terms of a semidefinite
programming (SDP) constraint. In what follows, we briefly describe the
basic ideas behind SDP and sum of squares programming, and their
applications to our problem.

\paragraph{Semidefinite programming.} SDP is a specific kind of convex
optimization problem with very appealing numerical properties. An SDP
problem corresponds to the optimization of a linear function over the
intersection of an affine subspace and the cone of positive
semidefinite matrices. For much more information about SDP and its
many applications, we refer the reader to the surveys
\cite{VaB:96,ToddSDP} and the comprehensive treatment in
\cite{HandSDP}.

An SDP problem in standard primal form is usually written as:
\begin{align*}
\mathrm{minimize}   \quad      C \bullet &X   \quad   &
\mbox{subject to}   \quad      A_i \bullet X  &= b_i, \quad i = 1,\ldots,m \\
                &           &          X &\succeq 0,
\end{align*}
where $C, A_i$ are symmetric $n \times n$ matrices, and $X \bullet Y
:= \mathrm{trace}(X Y)$. The symmetric matrix $X$ is the optimization
variable over which the maximization is performed.  The inequality in
the second line means that the matrix $X$ must be positive
semidefinite, i.e., all its eigenvalues should be greater than or
equal to zero.  The set of feasible solutions, i.e., the set of
matrices $X$ that satisfy the constraints, is always a convex set. In
the particular case when $C=0$, the problem reduces to whether or not
the inequality can be satisfied for some matrix $X$. In this case, the
SDP is referred to as a \emph{feasibility problem}.

There are a number of sophisticated and reliable methods to
numerically solve semidefinite programming problems. One of the most
successful approaches is based on \emph{primal-dual interior point
methods}, that generalize many of the techniques used in linear
programming \cite{NN}. The interior-point approach to SDP typically
involves the iterative solution of a perturbed version of the KKT
optimality conditions. Each iteration requires the computation of the
corresponding Newton direction, and the solution of a system of linear
equations. A theoretical bound on the number of Newton iterations is
$O(\sqrt{n} \log \frac{1}{\epsilon})$ for an $\epsilon$-approximate
solution. This estimate is signficantly more conservative than what is
usually experienced in practice, where the dependence on $n$ is very
mild (typically, 10-40 Newton iterations are enough for most
problems).  The cost of each iteration heavily depends on the
structure and sparsity of the matrices $A_i$, and is dominated by the
computation of the Hessian and the solution of the corresponding
linear system. In the fully dense case, this cost is of the order of
$\max\{mn^3,m^2n^2,m^3\}$, where the first two terms correspond to the
construction of the Hessian, and the last one to the solution of the
Newton system.

\paragraph{Sums of squares programming.}
Consider a given multivariate polynomial for which we want to decide
whether a sum of squares decomposition exists. This question is
equivalent to a semidefinite programming (SDP) problem, because of the
following result, that has appeared in different forms in the work of
Shor \cite{Shor}, Choi-Lam-Reznick \cite{ChoiLamReznick}, Nesterov
\cite{NesterovSquared}, and Parrilo \cite{Phd:Parrilo,sdprelax}.
\begin{theorem}
A homogeneous multivariate polynomial $p(x)$ of degree $2d$ is a sum
of squares if and only if
\begin{equation}
p(x) = (x^{[d]})^T Q x^{[d]},
\label{Par:sosrep}
\end{equation}
where $x^{[d]}$ is a vector whose entries are (possibly scaled)
monomials of degree $d$ in the variables $x_i$, and $Q$ is a symmetric
positive semidefinite matrix.
\end{theorem}
Since in general the entries of $x^{[d]}$ are not algebraically
independent, the matrix $Q$ in the representation (\ref{Par:sosrep})
\emph{is not unique}. In fact, there is an affine subspace of matrices
$Q$ that satisfy the equality, as can be easily seen by expanding the
right-hand side and equating term by term. To obtain an SOS
representation, we need to find a positive semidefinite matrix in this
affine subspace. Therefore, the problem of checking if a polynomial
can be decomposed as a sum of squares is \emph{equivalent} to
verifying whether a certain affine matrix subspace intersects the cone
of positive definite matrices, and hence an SDP feasibility problem.

\begin{example}
  Consider the quartic homogeneous polynomial in two variables
  described below, and define the vector of monomials as $[ x^2, y^2,
  x y]^T$.
\begin{eqnarray*}
p(x,y) &=& 2 x^4 + 2 x^3 y  - x^2 y^2 + 5 y^4 \\
&=&
\left[\begin{array}{c}
x^2 \\  y^2 \\ x y
\end{array}\right]^T
\left[\begin{array}{ccc}
q_{11} & q_{12} & q_{13} \\
q_{12} & q_{22} & q_{23} \\
q_{13} & q_{23} & q_{33}
\end{array}\right]
\left[\begin{array}{c}
x^2 \\  y^2 \\ x y
\end{array}\right]\\
&=&
q_{11} x^4 + q_{22} y^4 + (q_{33} + 2 q_{12}) x^2 y^2 + 2 q_{13} x^3 y + 2 q_{23} x y^3
\end{eqnarray*}
For the left- and right-hand sides to be identical, the following
linear equations should hold:
\begin{equation}
q_{11} = 2, \quad
q_{22} = 5, \quad
q_{33} + 2 q_{12} = -1, \quad
2 q_{13} = 2, \quad
2 q_{23} = 0.
\end{equation}

A positive semidefinite $Q$ that satisfies the linear equalities can
then be found using SDP. A particular solution is given by:
\[
Q =
\left[\begin{array}{rrr}
2  & -3 & 1 \\ -3 & 5 & 0 \\ 1 & 0 & 5
\end{array}\right]
= L^T L, \qquad
L =
\frac{1}{\sqrt{2}}\left[\begin{array}{rrr}
2 & -3 & 1 \\
0 & 1 & 3
\end{array}\right],
\]
and therefore we have the sum of squares decomposition:
\[
p(x,y) = \frac{1}{2} (2 x^2 - 3 y^2 + x y)^2 +
\frac{1}{2}(y^2 + 3 x y)^2.
\]
\label{Par:ex:sosexample}
\hfill $\square$
\end{example}

\subsection{Norms and SOS polynomials}

The procedure described in the previous subsection can be easily
adapted to the case where the polynomial $p(x)$ is not fixed, but
instead we search for an SOS polynomial in a given affine family (for
instance, all homogeneous polynomials of a given degree).

This line of thought immediately suggests the following SOS relaxation
of the conditions in Theorem~\ref{thm:psdbound}:
\begin{equation}
\rho_{SOS,2d} :=
\inf_{p(x) \in \R_{2d}[x], \gamma} \gamma \qquad \mbox{s.t. }\left\{
\begin{array}{rl}
p(x) \, & \mbox{is SOS}\\
\gamma^{2d} p(x) - p(A_i x) \, & \mbox{is SOS}
\end{array}
\right.
\label{eq:SOSrelax}
\end{equation}
where $\R_{2d}[x]$ is the set of homogeneous polynomials of degree
$2d$.

\begin{remark}
Theorem~\ref{thm:psdbound} requires a strictly positive polynomial
$p(x)$, so it would be natural to add some strict positivity condition
to the relaxation~(\ref{eq:SOSrelax}). For instance, one could require
for the polynomial $p(x)$ to belong to the relative interior of the
SOS cone.  However, since interior-point methods by construction
always produce solutions in the relative interior of the corresponding
convex set, this is automatically satisfied if the problem is
feasible. Alternatively, it is possible to give a formulation that
includes terms of the form $\epsilon ||x||^{2d}$, for small positive
$\epsilon$. These modifications are unnecessary in practice.
\end{remark}

For any fixed degree $d$ and any given $\gamma$, the constraints in
this problem are all of SOS type, and thus equivalent to semidefinite
programming. Therefore, the computation of $\rho_{SOS,2d}$ is a
quasiconvex problem, and can be easily solved with a standard SDP
solver, and a simple bisection method for the scalar variable
$\gamma$. By Theorem~\ref{thm:psdbound}, the solution of this
relaxation yields an upper bound on the joint spectral radius
\begin{equation}
\rho(A_1,\ldots,A_m) \leq \rho_{SOS,2d},
\label{eq:trivbound}
\end{equation}
where $2d$ is the degree of the approximating polynomial.

\subsection{Quality of approximation}

What can be said about the quality of the bounds produced by the SOS
relaxation? We present next some results to answer this question; a
more complete characterization is developed in
Section~\ref{sec:goodbounds}. An inspiring result in this direction is
the following theorem of Barvinok, that quantifies how tightly SOS
polynomials can approximate norms:
\begin{theorem}[\cite{Barvinok}, p.~221]
\label{thm:Barvinok}
Let $||\cdot||$ be a norm in $\R^n$. For any integer $d \geq 1$ there
exists a homogeneous polynomial $p(x)$ in $n$ variables of degree $2d$
such that
\begin{enumerate}
\item The polynomial $p(x)$ is a sum of squares.
\item For all $x \in \R^n$,
\[
p(x)^\frac{1}{2d} \leq ||x|| \leq k(n,d) \,
p(x)^\frac{1}{2d},
\]
where $k(n,d) := \binom{n+d-1}{d}^{\frac{1}{2d}}$.
\end{enumerate}
\end{theorem}
For fixed state dimension $n$, by increasing the degree $d$ of the
approximating polynomials, the factor in the upper bound can be made
arbitrarily close to one. In fact, for large $d$, we have the
approximation
\[
k(n,d) \; \approx \; 1 + \frac{n-1}{2} \frac{\log d}{d}.
\]

To apply these results to our problem, consider the following. If
$\rho(A_1,\ldots,A_m) < \gamma$, by Theorem~\ref{thm:RotaStrang} (and
sharper results in \cite{Bar88,Koz90,wirth}) there exists a norm
$\|\cdot\|$ such that
\[
||A_i x|| \leq \gamma ||x||, \quad \forall x \in \R^n, i = 1,\ldots,m.
\]
By Theorem~\ref{thm:Barvinok}, we can therefore approximate this norm
with a homogeneous SOS polynomial $p(x)$ of degree $2d$ that will then
satisfy
\[
p(A_i x)^\frac{1}{2d}\leq ||A_i x|| \leq \gamma ||x||
\leq \gamma \, k(n,d) \, p(x)^\frac{1}{2d},
\]
and thus we know that there exists a feasible solution of
\[
\left\{
\begin{array}{rl}
p(x) \, & \mbox{is SOS}\\
\alpha^{2d} p(x) - p(A_i x) \, & \geq 0 \qquad i=1,\ldots,m,
\end{array}
\right.
\]
for $\alpha = k(n,d) \rho(A_1,\ldots,A_m)$.

Despite these appealing results, notice that in general we cannot yet
conclude from this that the proposed SOS relaxation will always obtain
a solution that is within $k(n,d)^{-1}$ from the true spectral
radius. The reason is that even though we can prove the existence of a
$p(x)$ that is SOS and for which $\alpha^{2d} p(x) - p(A_i x)$ are
nonnegative for all $i$, it is unclear whether the last $m$
expressions are actually SOS. We will show later in the paper that
this is indeed the case. Before doing this, we concentrate first on
two important cases of interest, where the described approach
guarantees a good quality of approximation.

\paragraph{Planar systems.}
The first case corresponds to two-dimensional (planar) systems, i.e.,
when $n=2$. In this case, it always holds that nonnegative homogeneous
bivariate polynomials are SOS (e.g., \cite{Reznick}). Thus, we have
the following result:
\begin{theorem}
Let $\{A_1,\ldots,A_m\} \subset \R^{2 \times 2}$. Then, the SOS
relaxation~(\ref{eq:SOSrelax}) always produces a solution satisfying:
\[
{\textstyle \frac{1}{2}} \rho_{SOS,2d} \leq
(d+1)^{-\frac{1}{2d}} \,
\rho_{SOS,2d} \leq \rho(A_1,\ldots,A_m) \leq \rho_{SOS,2d}.
\]
This result is \emph{independent} of the number $m$ of matrices.
\end{theorem}

\paragraph{Quadratic Lyapunov functions.}
In the quadratic case (i.e., $2d=2$), it is also true that nonnegative
quadratic forms are sums of squares. Since
\[
{\binom{n+d-1}{d}}^\frac{1}{2d} =
\binom{n}{1}^\frac{1}{2} = \sqrt{n},
\]
the inequality
\begin{equation}
\frac{1}{\sqrt{n}} \; \rho_{SOS,2} \leq \rho(A_1,\ldots,A_m) \leq \rho_{SOS,2}
\label{eq:quadlyapbound}
\end{equation}
follows. This bound exactly coincides with the results of Ando and
Shih \cite{Ando98} or Blondel, Nesterov and Theys \cite{BlNT04}. This
is perhaps not surprising, since in this case both Ando and Shih's
proof \cite{Ando98} and Barvinok's theorem rely on the use of John's
ellipsoid to approximate the same underlying convex set.

\paragraph{Level sets and convexity}
Unlike the norms that appear in Theorem~\ref{thm:RotaStrang}, an
appealing feature of the SOS-based method is that we are not
constrained to use polynomials with convex level sets. This enables in
some cases much better bounds than what is promised by the theorems
above, as illustrated in the following example.

\begin{figure}[t]
\centering
\includegraphics[width=0.5\columnwidth]{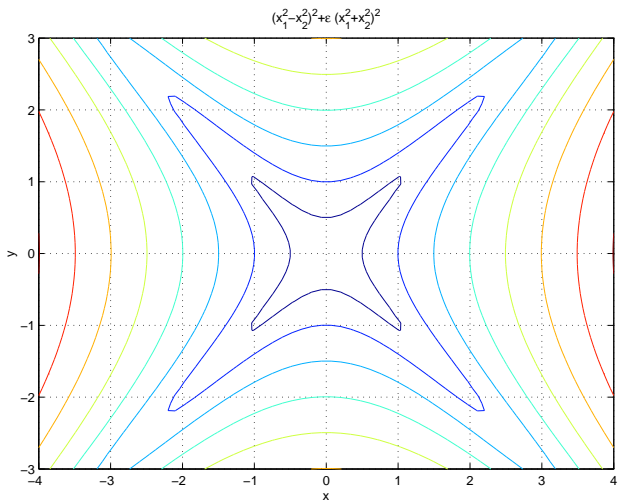}
\caption{Level sets of the quartic homogeneous polynomial
$V(x_1,x_2)$. These define a Lyapunov function, under which both $A_1$
and $A_2$ are $(1+\epsilon)$-contractive. The value of $\epsilon$ is
here equal to $0.01$.}
\label{fig:jsr}
\end{figure}

\begin{example}
This is based on a construction by Ando and Shih
\cite{Ando98}. Consider the problem of proving a bound on the joint
spectral radius of the following matrices:
\[
A_1 =
\left[\begin{array}{cc}
1 & 0 \\ 1 & 0
\end{array}\right], \qquad
A_2 =
\left[\begin{array}{rr}
0 & 1 \\ 0 & -1
\end{array}\right].
\]
For these matrices, it can be easily shown that
$\rho(A_1,A_2)=1$. Using a common quadratic Lyapunov function (i.e.,
the case $d=2$), the upper bound on the joint spectral radius is equal
to $\sqrt{2}$. However, a simple quartic SOS Lyapunov function is
enough to prove an upper bound of $1+\epsilon$ for every $\epsilon
>0$, since the SOS polynomial
\[
V(x) = (x_1^2-x_2^2)^2 + \epsilon (x_1^2+x_2^2)^2
\]
satisfies
\begin{eqnarray*}
(1+\epsilon) V(x) - V(A_1 x) &=& ( x_2^2-x_1^2+\epsilon (x_1^2+x_2^2) )^2 \\
(1+\epsilon) V(x) - V(A_2 x) &=& ( x_1^2-x_2^2+\epsilon (x_1^2+x_2^2) )^2.
\end{eqnarray*}
The corresponding level sets of $V(x)$ are plotted in
Figure~\ref{fig:jsr}, and are clearly non-convex.
\label{ex:ando}
\end{example}

\section{Symmetric algebra and induced matrices}
\label{sec:symmalgebra}

We present next some further bounds on the quality of the SOS
relaxation~(\ref{eq:SOSrelax}), either by a more refined analysis of
the SOS polynomials in Barvinok's theorem or by explicitly producing
an SOS Lyapunov function of guaranteed suboptimality properties. These
constructions are quite natural, and parallel some lifting ideas as
well as the classical iteration used in the solution of discrete-time
Lyapunov inequalities. Before proceeding further, we briefly revisit
some classical notions from multilinear algebra.

\paragraph{Symmetric algebra of a vector space}
Consider a vector $x \in \R^n$, and an integer $d \geq 1$. We define
its $d$-lift $x^{[d]}$ as a vector in $\R^N$, where $N: =
\binom{n+d-1}{d}$, with components $\{ \sqrt{\alpha !} \, x^\alpha
\}_\alpha$, where $\alpha = (\alpha_1,\ldots,\alpha_n)$, $|\alpha| :=
\sum_i \alpha_i = d$, and $\alpha !$ denotes the multinomial
coefficient $\alpha ! := \binom{d}{\alpha_1,\alpha_2,\ldots,\alpha_n}=
\frac{d!}{\alpha_1!  \alpha_2! \ldots \alpha_n!}$. That is, the
components of the lifted vector are the monomials of degree $d$,
scaled by the square root of the corresponding multinomial
coefficients.
\begin{example}
Let $n=2$, and $x = [u,v]^T$. Then, we have
\[
\left[\begin{array}{c} u \\ v \end{array}\right]^{[1]} =
\left[\begin{array}{c} u \\ v \end{array}\right], \qquad
\left[\begin{array}{c} u \\ v \end{array}\right]^{[2]} =
\left[\begin{array}{c} u^2 \\ \sqrt{2} u v \\ v^2 \end{array}\right], \qquad
\left[\begin{array}{c} u \\ v \end{array}\right]^{[3]} =
\left[\begin{array}{c} u^3 \\ \sqrt{3} u^2 v \\
\sqrt{3} u v^2 \\ v^3 \end{array}\right].
\]
\end{example}
The main motivation for this specific scaling of the components, is to
ensure that the lifting preserves some of the properties of the
underlying normed space. In particular, if $||\cdot||$ denotes the
standard Euclidean norm, it can be easily verified that $||x^{[d]}|| =
||x||^d$. Thus, the lifting operation provides a norm-preserving (up
to power) embedding of $\R^n$ into $\R^N$. When the original space is
projective, this is the so-called \emph{Veronese} embedding.

This concept can be directly extended from vectors to linear
transformations. Consider a linear map in $\R^n$, and the associated
$n \times n$ matrix $A$. Then, the lifting described above naturally
induces an associated map in $\R^N$, that makes the corresponding
diagram commute.  The matrix representing this linear transformation
is the \emph{$d$-th induced matrix} of $A$, denoted by $A^{[d]}$,
which is the unique $N \times N$ matrix that satisfies
\[
A^{[d]} x^{[d]} = (A x)^{[d]}.
\]
In systems and control, these classical constructions of multilinear
algebra have been used under different names in several works, among
them \cite{BrockettLie,Zelen} and (implicitly) \cite{BlNes05}.
Although not mentioned in the Control literature, there exists a
simple explicit formula for the entries of these induced matrices; see
\cite{MarcusMultilinear,MarcusMinc}.  The $d$-th induced matrix
$A^{[d]}$ has dimensions $N \times N$. Its entries are given by
\begin{equation}
(A^{[d]})_{\alpha \beta} = \frac{\mathrm{per}\,  A(\alpha,\beta)}{\sqrt{\mu(\alpha) \mu(\beta)}},
\label{eq:perm}
\end{equation}
where the indices $\alpha,\beta$ are all the $d$-element multisets of
$\{1,\ldots,n\}$, the notation $\mathrm{per}$ indicates the
\emph{permanent}\footnote{The permanent of a matrix $A \in \R^{n
\times n}$ is defined as $\textrm{per}(A):=\sum_{\sigma \in \Pi_n}
\prod_{i=1}^n a_{i,\sigma(i)}$, where $\Pi_n$ is the set of all
permutations in $n$ elements.} of a square matrix, and $\mu(S)$ is the
product of the factorials of the multiplicities of the elements of the
multiset $S$.
\begin{example}
Consider the case $n=2$, $d=3$. The corresponding 3-element multisets
are $\{1,1,1\}$, $\{1,1,2\}$, $\{1,2,2\}$ and $\{2,2,2\}$. The third
induced matrix is then
\begin{align*}
A^{[3]} &=
\begin{bmatrix}
                       a_{11}^3&          \sqrt{3} a_{11}^2 a_{12} &          \sqrt{3} a_{11} a_{12}^2 &                       a_{12}^3 \\
          \sqrt{3} a_{11}^2 a_{21}& a_{11} (a_{11} a_{22}+2 a_{21} a_{12})& a_{12} (2 a_{11} a_{22}+a_{21} a_{12}) &          \sqrt{3} a_{12}^2 a_{22} \\
          \sqrt{3} a_{11} a_{21}^2& a_{21} (2 a_{11} a_{22}+a_{21} a_{12})& a_{22} (a_{11} a_{22}+2 a_{21} a_{12}) &          \sqrt{3} a_{12} a_{22}^2 \\
                       a_{21}^3&          \sqrt{3} a_{21}^2 a_{22}  &          \sqrt{3} a_{21} a_{22}^2 &                       a_{22}^3 \\
\end{bmatrix}.
\end{align*}
\end{example}
It can be shown that these operations define an algebra homomorphism,
i.e., they respect the structure of matrix multiplication. In
particular, for any matrices $A,B$ of compatible dimensions, the
following identities hold:
\[
(A B)^{[d]} = A^{[d]} B^{[d]}, \qquad (A^{-1})^{[d]} = (A^{[d]})^{-1}.
\]
Furthermore, there is a simple and appealing relationship between the
eigenvalues of $A^{[d]}$ and those of $A$. Concretely, if
$\lambda_1,\ldots,\lambda_n$ are the eigenvalues of $A$, then the
eigenvalues of $A^{[d]}$ are given by $\prod_{j \in S} \lambda_j$
where $S \subseteq \{1,\ldots,n\}, |S| = d$; there are exactly
$\binom{n+d-1}{d}$ such multisets. A similar relationship holds for
the corresponding eigenvectors. Essentially, as explained below in
more detail, the induced matrices are the symmetry-reduced
version of the $d$-fold Kronecker product.

The symmetric algebra and associated induced matrices are classical
objects of multilinear algebra. Induced matrices, as defined above, as
well as the more usual \emph{compound matrices}, correspond to two
specific isotypic components of the decomposition of the $d$-fold
tensor product under the action of the symmetric group $S^d$ (i.e.,
the \emph{symmetric} and \emph{skew-symmetric} algebras).  Compound
matrices are associated with the alternating character (hence their
relationship with determinants), while induced matrices correspond
instead to the trivial character, thus the connection with
permanents. Similar constructions can be given for any other character
of the symmetric group, by replacing the permanent in (\ref{eq:perm})
with the suitable immanants; see \cite{MarcusMultilinear} for
additional details.

\subsection{Bounds on the quality of $\rho_{SOS,2d}$}
\label{sec:goodbounds}

In this section we present a bound on the approximation properties of
the SOS approximation, based on the ideas introduced above. As we will
see, the techniques based on the lifting described will exactly yield
the factor $k(n,d)^{-1}$ suggested by Barvinok's theorem.

We first prove a preliminary result on the behavior of the joint
spectral radius under $d$-lifting.  The scaling properties described
earlier can be applied to obtain the following:
\begin{lemma}
Given matrices $\{A_1,\ldots,A_m\} \subset \R^{n \times n}$ and an
integer $d \geq 1$, the following identity holds:
\[
\rho(A_1^{[d]},\ldots,A_m^{[d]}) = \rho(A_1,\ldots,A_m)^d.
\]
\label{lem:scalingjsr}
\end{lemma}
The proof follows directly from the definition~(\ref{eq:defjsr}) and
the two properties $(A B)^{[d]} = A^{[d]} B^{[d]}$, $||x^{[d]}|| =
||x||^d$, and it is thus omitted.

Combining all these inequalities, we obtain the main result of this paper:
\begin{theorem}
The SOS relaxation (\ref{eq:SOSrelax}) satisfies:
\begin{equation}
{\textstyle\binom{n+d-1}{d}}^{-\frac{1}{2d}} \; \rho_{SOS,2d} \leq \rho(A_1,\ldots,A_m) \leq \rho_{SOS,2d}.
\label{eq:sos2dbound}
\end{equation}
\label{thm:sos2dbound}
\end{theorem}
\begin{proof}
Since the dimension of $A_i^{[d]}$ is $\binom{n+d-1}{d}$, from
Lemma~\ref{lem:scalingjsr} and inequality (\ref{eq:quadlyapbound}) it
follows that:
\[
{\textstyle\binom{n+d-1}{d}}^{-\frac{1}{2}}
\; \rho_{SOS,2}(A_1^{[d]},\ldots,A_m^{[d]})
\leq
\rho(A_1^{[d]},\ldots,A_m^{[d]}) = \rho(A_1,\ldots,A_m)^d.
\]
Combining this with (\ref{eq:trivbound}) and the inequality (proven
later in Theorem~\ref{thm:3bounds}),
\[
\rho_{SOS,2d}(A_1,\ldots,A_m)^d \leq \rho_{SOS,2}(A_1^{[d]},\ldots,A_m^{[d]}),
\]
the result follows.
\end{proof}

\section{Sum of squares Lyapunov iteration}
\label{sec:soslyap}

We describe next an alternative approach to obtain bounds on the
quality of the SOS approximation. As opposed to the results in the
previous section, the bounds now explicitly depend on the number of
matrices, but will usually be tighter in the case of small $m$.

Consider the iteration defined by
\begin{equation}
V_0(x) = 0, \qquad V_{k+1}(x) = Q(x) + \frac{1}{\beta} \sum_{i=1}^m V_k(A_i x),
\label{eq:iteration}
\end{equation}
where $Q(x)$ is a fixed $n$-variate homogeneous polynomial of degree
$2d$ and $\beta > 0$.  The iteration defines an affine map in the
space of homogeneous polynomials of degree $2d$. As usual, the
iteration will converge under certain assumptions on the spectral
radius of this linear operator.
\begin{theorem}
The iteration defined in (\ref{eq:iteration}) converges for arbitrary
$Q(x)$ if $\rho(A_1^{[2d]} + \cdots + A_m^{[2d]}) < {\beta}$.
\label{thm:convergence}
\end{theorem}
\begin{proof}
The vector space of homogenous polynomials $\R_{2d}[x_1,\ldots,x_n]$
is naturally isomorphic to the space of linear functionals on
$(\R^n)^{[2d]}$, via the identification $V_k(x) = \langle v_k ,
x^{[2d]} \rangle$, where $v_k \in \R^{\binom{n+2d-1}{2d}}$ is the
vector of (scaled) coefficients of $V_k(x)$. Then, since $V_k(A_i x) =
\langle v_k, (A_i x)^{[2d]} \rangle = \langle v_k, A_i^{[2d]} x^{[2d]}\rangle=
\langle (A_i^{[2d]})^T v_k, x^{[2d]}\rangle$, the iteration
(\ref{eq:iteration}) can be simply expressed as:
\[
v_{k+1} = q + \frac{1}{\beta} \left( \sum_{i=1}^m A_i^{[2d]} \right)^T v_{k},
\]
and it is well known that an affine iteration converges if the
spectral radius of the linear term is less than one.
\end{proof}

For simplicity of notation, we define the following quantity,
corresponding to the spectral radius of the sum of the $2d$-lifted
matrices:
\begin{equation}
\rho_{SR,2d} := \rho(A_1^{[2d]} + \cdots + A_m^{[2d]})^\frac{1}{2d}.
\label{eq:rhold}
\end{equation}

\begin{theorem}
\label{thm:sosvsnesterov}
The following inequality holds:
\[
\rho_{SOS,2d} \leq \rho_{SR,2d}
\]
\end{theorem}
\begin{proof}
Choose a $Q(x)$ that is in the interior of the SOS cone, e.g., $Q(x)
:= (\sum_{i=1}^n x_i^2)^d$, and let $\beta = \rho(A_1^{[2d]} + \cdots
+ A_m^{[2d]})+\epsilon$. The iteration~(\ref{eq:iteration}) guarantees
that $V_{k+1}$ is SOS if $V_{k}$ is. By induction, all the iterates
$V_k$ are SOS.  By the choice of $\beta$ and
Theorem~\ref{thm:convergence}, the $V_k$ converge to some homogeneous
polynomial $V_\infty(x)$. By the closedness of the cone of SOS
polynomials, the limit $V_\infty$ is also SOS.  Furthermore, we have
\[
\beta V_\infty(x) - V_\infty(A_i x) = \beta Q(x) + \sum_{j \not = i} V_\infty (A_j x)
\]
and therefore the expression on the left-hand side is SOS. This
implies that $p(x):=V_\infty(x)$ is a feasible solution of the SOS
relaxation (\ref{eq:SOSrelax}). Taking $\epsilon \rightarrow 0$, the
result follows.
\end{proof}
Notice that if the spectral radius condition in
Theorem~\ref{thm:convergence} is satisfied, then for any fixed $Q(x)$
the corresponding limit $V_\infty(x) = \langle v_\infty,
x^{[2d]}\rangle$ can be simply obtained by solving the nonsingular
system of linear equations
\[
\left(I- \frac{1}{\beta}\sum_{i=1}^m A_i^{[2d]} \right)^T v_\infty = q,
\]
thus generalizing the standard Lyapunov equation. The iteration
argument is only used to prove that the solution of this linear system
yields a strictly positive SOS polynomial. A slightly different
approach here is via the finite-dimensional version of the
Krein-Rutman theorem (or generalized Perron-Frobenius); see for
instance \cite{Protasov1} or \cite{ParriloKhatri}.

\begin{theorem}
The SOS relaxation (\ref{eq:SOSrelax}) satisfies:
\[
m^{-\frac{1}{2d}} \, \rho_{SOS,2d} \leq \rho(A_1,\ldots,A_m) \leq \rho_{SOS,2d}.
\]
\label{thm:msos2dbound}
\end{theorem}
\begin{proof}
This follows directly from inequality~(\ref{eq:trivbound}), and the fact that
\[
\rho_{SOS,2d} \leq \rho\left(\sum_{i=1}^m A_i^{[2d]}\right)^\frac{1}{2d} \\
\leq  m^\frac{1}{2d} \cdot \rho (A_1^{[2d]}, \ldots,  A_m^{[2d]})^\frac{1}{2d} \\
= m^\frac{1}{2d} \cdot \rho \left(A_1, \ldots , A_m \right),
\]
where the first inequality is Theorem~\ref{thm:sosvsnesterov}, the
second one follows from the general fact that $\rho(A_1+\cdots+A_m)
\leq m \rho(A_1,\ldots,A_m)$ (see e.g., Corollary 1 in \cite{BlNes05}), and
the third from Lemma~\ref{lem:scalingjsr}.
\end{proof}
The iteration~(\ref{eq:iteration}) is the natural generalization of
the Lyapunov recursion for the single matrix case, and of the
construction by Ando and Shih in \cite{Ando98} for the quadratic
case. By the remarks in Section~\ref{sec:symmalgebra} above, and as
described in more detail in the next section, it can be shown that the
quantity $\rho_{SR,2d}$ is essentially equal to those defined by
Protasov in \cite[\S 4]{Protasov1} and Blondel and Nesterov in
\cite{BlNes05}. As a consequence of Theorem~\ref{thm:sosvsnesterov},
the SOS-based approach will \emph{always} produce estimates at least
as good as the ones given by these procedures.

\section{Comparison with earlier techniques}
\label{sec:comparison}

In this section we compare the $\rho_{SOS,2d}$ approach with some
earlier bounds from the literature. We show that our bound is never
weaker than those obtained by all the other procedures.

\subsection{Methods of Protasov and Blondel-Nesterov}

Protasov \cite{Protasov1} has shown that an upper bound on the
``standard'' joint spectral radius can be computed via the so-called
joint $p$-radius, a generalization of the definition~(\ref{eq:defjsr})
involving $p$-norms. Furthermore, he has shown that in the case of
even integer $p$, the value of the $p$-radius of an irreducible finite
set of matrices exactly corresponds to the spectral radius of a single
operator, that can in principle be constructed based on the matrices
$A_i$.

Independently, Blondel and Nesterov \cite{BlNes05} developed a
technique based on the calculation of the spectral radius of
``lifted'' matrices. In fact, they present two different lifting
procedures (``Kronecker'' and ``semidefinite'' liftings), and in
Section~5 of their paper, they describe a family of bounds obtained by
arbitrary combinations of these two liftings.

Both of these methods are in fact equivalent to our construction of
$\rho_{SR,2d}$ in Section~\ref{sec:soslyap}, in the sense that they
all yield exactly the same numerical value. By
Theorem~\ref{thm:sosvsnesterov}, they are thus also weaker than the
SOS-based construction.  The bound defined by $\rho_{SR,2d}$
in~(\ref{eq:rhold}) relies on a single canonically defined lifting,
and requires much less numerical effort than the Blondel-Nesterov
construction. Furthermore, instead of the somewhat more complicated
construction of Protasov, the expression of the entries of the lifted
matrices are given by the simple formula~(\ref{eq:perm}), making a
computer implementation straightforward, with no irreducibility
assumptions being required.

It can be shown that our construction (or Protasov's) exactly
corresponds to a fully symmetry-reduced version of the
Blondel-Nesterov procedure, thus yielding equivalent bounds, but at a
much smaller computational cost since the corresponding matrices are
exponentially smaller (for fixed $n$, the size grows as $O(d^{n-1})$
as opposed to $O(n^{2d})$). Therefore, even if no SDPs are to be
solved (as would be required by the tighter bound $\rho_{SOS,2d}$),
the formulation in terms of the matrices $A_i^{[2d]}$ still has many
advantages.

\begin{table}[t]
\begin{center}
\begin{tabular}{|c|c || c|c || c|c || c|c|}
\hline
& & \multicolumn{2}{c||}{ \cite{BlNes05}, Kronecker } & 
\multicolumn{2}{c||}{ \cite{BlNes05}, semidefinite } & 
\multicolumn{2}{c|}{This paper} \\
\cline{3-8}
Steps / $2d$ &Accuracy&$n=2$& $n=10$ &$n=2$ & $n=10$ & $n=2$ & $n=10$ \\
\hline\hline
1 / 2  & 0.707 & 4 & 100 & 3 & 55 & 3 & 55  \\
2 / 4  & 0.840 & 16 & 10000 & 6 & 1540 & 5 & 715 \\
3 / 8  & 0.917 & 256 & $10^8$ & 21 & 1186570 & 9 & 24310 \\
4 / 16 & 0.957 & 65536 & $10^{16}$ & 231 & $7.04 \times 10^{11}$ & 17 & 2042975 \\
5 / 32 & 0.978 & $4.29\times 10^9$ & $10^{32}$ &   26796 & $2.48 \times 10^{23}$& 33 & $3.5 \times 10^8$  \\
\hline
\end{tabular}
\end{center}
\caption{Comparison of matrix sizes for the different lifting
procedures to compute $\rho_{SR,2d}$. The matrix size for the
Kronecker lifting is $n^{2d}$, while the recursive semidefinite
lifting is given by the $d$-step recursion $s_{2k} = \binom{s_k+1}{2}$
with $s_1=n$, and the size for the symmetric algebra approach is
$\binom{n+2d-1}{2d}$. The accuracy estimates correspond to the case of
two matrices, i.e., $m=2$.}
\label{tab:BNtwo}
\end{table}
As an illustrative comparison of the advantages of this reduced
formulation, in Table~\ref{tab:BNtwo} we present the sizes of the
matrices required by the method in~\cite{BlNes05} (using the
``Kronecker'' and ``recursive semidefinite'' liftings) and our
approach to $\rho_{SR,2d}$ via the symmetric algebra. The data in
Table~\ref{tab:BNtwo} corresponds to that in~\cite[p.~266]{BlNes05}
(with a minor misprint corrected).

\subsection{Common quadratic Lyapunov functions}

This method corresponds to finding a common quadratic Lyapunov
function, either directly for the matrices $A_i$, or for the lifted
matrices $A_i^{[d]}$. Specifically, let
\[
\rho_{CQ,2d} := \inf \, \left \{ \; \gamma \; \; | \; \; \gamma^{2d} P -
(A_i^{[d]})^T P A_i^{[d]} \succeq 0, \quad P \succ 0 \right \}.
\]
This is essentially equivalent to what is discussed in Corollary 3 of
\cite{BlNes05}, except that the matrices involved in our approach are
exponentially smaller (of size $\binom{n+d-1}{d}$ rather than $n^d$),
as all the symmetries have been taken out\footnote{There seems to be a
typo in equation (7.4) of \cite{BlNes05}, as all the terms $A_i^k$
should likely read $A_i^{\otimes k}$.}. Notice also that, as a
consequence of their definitions, we have
\[
\rho_{CQ,2d}(A_1,\ldots,A_m)^d = \rho_{SOS,2}(A_1^{[d]},\ldots,A_m^{[d]}).
\]

We can then collect most of these results in a single theorem:
\begin{theorem}
The following inequalities between all the bounds hold:
\begin{equation}
\rho(A_1,\ldots,A_m) \leq \rho_{SOS,2d} \leq \rho_{CQ,2d} \leq
\rho_{SR,2d}.
\label{eq:ineqs}
\end{equation}
\label{thm:3bounds}
\end{theorem}
\begin{proof}
The left-most inequality is~(\ref{eq:trivbound}). The right-most
inequality follows from a similar (but stronger) argument to the one
given in Theorem~\ref{thm:sosvsnesterov} above, since the spectral
radius condition $\rho(A_1^{[2d]}+\cdots + A_m^{[2d]})< \beta$
actually implies the convergence of the matrix iteration in
$\mathcal{S}^{N}$ given by
\[
P_{k+1} = Q + \frac{1}{\beta} \sum_{i=1}^m (A_i^{[d]})^T P_k A_i^{[d]}, \qquad P_0 = I.
\]

For the middle inequality, let $p(x):= (x^{[d]})^T P
x^{[d]}$. Since $P \succ 0$, it follows that $p(x)$ is SOS. From
$\gamma^{2d} P - (A_i^{[d]})^T P A_i^{[d]} \succeq 0$, left- and
right-multiplying by $x^{[d]}$, we have that $\gamma^{2d} p(x) - p(A_i
x)$ is also SOS, and thus $p(x)$ is a feasible solution
of~(\ref{eq:SOSrelax}), from where the result directly follows.
\end{proof}

\begin{remark}
We always have $\rho_{SOS,2} = \rho_{CQ,2}$, since both correspond
to the case of a common quadratic Lyapunov function for the matrices $A_i$.
\end{remark}

\subsection{Computational cost}
In this section we quantify the computational cost of the bound
$\rho_{SOS,2d}$. In the following calculations we keep $d$ fixed, and
study the scaling behavior as a function of the dimension $n$.

As mentioned in Section~\ref{sec:sosnorms}, solving a semidefinite
programming problem typically requires several Newton iterations, with
the cost of each iteration being dominated by the construction of the
Hessian and solution of the corresponding linear system. For the SOS
bound $\rho_{SOS,2d}$, the underlying SDP problem has $m+1$ matrix
inequalities corresponding to the SOS constraints
in~(\ref{eq:SOSrelax}), each of dimension $\binom{n+d-1}{d} \approx
\frac{1}{d!} \cdot n^d$, which is $O(n^d)$ for fixed $d$. The number of
decision variables is approximately $m \cdot \binom{n+2d-1}{2d}
\approx m \cdot n^{2d}$. Thus, using a simple bisection method for
$\gamma$, exploiting the block-diagonal structure, and the fact that
the number of Newton iterations is essentially constant, we obtain
that the approximate cost of obtaining an $\epsilon$-approximate
solution of $\rho_{SOS,2d}$ is $O(m \cdot n^{6d} \cdot \log
\frac{1}{\epsilon})$, where $d$ is chosen such that $\epsilon \approx
\frac{n}{2} \frac{\log d}{d}$ or $\epsilon \approx m^{-\frac{1}{2d}}$,
depending on whether we use bounds that depend on the number of
matrices (Theorem~\ref{thm:msos2dbound}) or not (Theorem
\ref{thm:sos2dbound}).

We remark that these quantities are a relatively coarse estimate of
the best possible algorithmic complexity, since very little structure
of the corresponding SDP problem is being exploited. It is known that
for structured problems such as the ones appearing here much more
efficient SDP-based algorithms can be developed. In particular, in the
context of sum of squares problems several techniques are known to
exploit some of the available structure for more efficient
computation; see \cite{GHNV,LofbergParrilo,RohVandenberghe}.

\subsection{Examples}
We present next two numerical examples that compare the described
techniques. In particular, we show that the bounds in
Theorem~\ref{thm:3bounds} can all be strict.

\begin{example}
Here we revisit the construction presented earlier in
Example~\ref{ex:ando}. For the matrices given there we have:
\begin{align*}
\rho_{SOS,2} &= \sqrt{2}, &
\rho_{CQ,2}&= \sqrt{2}, &
\rho_{SR,2d} &= \sqrt[2d]{2},
\\
\rho_{SOS,4} &= 1, &
\rho_{CQ,4}&= 1. &
\end{align*}
\end{example}

\begin{example}
\label{ex:threemats}
Consider the three $4 \times 4 $ matrices (randomly generated) given by:
\[
A_1 =
\left[
\begin{array}{rrrr}
     0  &   1  &    7  &    4 \\
     1  &   6  &   -2  &   -3 \\
    -1  &  -1  &   -2  &   -6 \\
     3  &   0  &    9  &    1
\end{array}
\right],
\quad
A_2 =
\left[
\begin{array}{rrrr}
    -3  &    3  &    0  &   -2\\
    -2  &    1  &    4  &    9\\
     4  &   -3  &    1  &    1\\
     1  &   -5  &   -1  &   -2
\end{array}
\right],
\quad
A_3 =
\left[
\begin{array}{rrrr}
     1   &   4  &    5  &   10 \\
     0   &   5  &    1  &   -4 \\
     0   &  -1  &    4  &    6 \\
    -1   &   5  &    0  &    1
\end{array}
\right].
\]
The value of the different approximations are presented in
Table~\ref{tab:comparison}. A lower bound is $\rho(A_1
A_3)^\frac{1}{2} \approx 8.9149$, which is extremely close (and
perhaps exactly equal) to the upper bound $\rho_{SOS,4}$. Notice from
the $d=2$ entry of Table~\ref{tab:comparison} that all the
inequalities~(\ref{eq:ineqs}) can be strict.

\begin{table}[t]
\begin{center}
\begin{tabular}{|c|cc|ccc|}
\hline $d$ & $\dim A_i^{[d]}$ & $\dim A_i^{[2d]}$ & $\rho_{SOS,2d}$ &
$\rho_{CQ,2d}$ & $\rho_{SR,2d}$ \\ \hline 1 & 4 & 10 & 9.761 & 9.761 &
12.519 \\ 2 & 10 & 35 & 8.92 & 9.01 & 9.887 \\ 3 & 20 & 84 & 8.92 &
8.92 & 9.3133 \\ \hline
\end{tabular}
\end{center}
\caption{Comparison of the different approximations for Example~\ref{ex:threemats}.}
\label{tab:comparison}
\end{table}
\end{example}

\section{Conclusions}
\label{sec:conclusions}

We introduced a novel scheme for the approximation of the joint
spectral radius of a set of matrices using sum of squares
programming. The method is based on the use of a multivariate
polynomial to provide a norm-like quantity under which all matrices
are contractive. We provided an asymptotically tight estimate for the
quality of the bound, which is independent of the number of
matrices. We also proposed an alternative bound, that depends on the
number $m$ of matrices, based on a generalization of a Lyapunov
iteration.

Our results can be alternatively interpreted in a simpler way as
providing a trajectory-preserving lifting to a higher dimensional
space, and proving contractiveness with respect to an ellipsoidal norm
in that space. In this case, a weaker estimate can be obtained by
computing the spectral radius of a fixed matrix.  These results
generalize earlier work of Ando and Shih~\cite{Ando98}, Blondel,
Nesterov and Theys~\cite{BlNT04}, and provide an improvement over the
lifting procedure of Blondel and Nesterov~\cite{BlNes05}. The good
performance of our procedure was also verified using numerical
examples.

\paragraph{Acknowledgement}
We thank the referees for their careful reading of the manuscript, and
their many useful suggestions.

\bibliographystyle{alpha}
\bibliography{jsr}

\begin{thebibliography}{GHND03}

\bibitem[AS98]{Ando98}
T.~Ando and M.-H. Shih.
\newblock Simultaneous contractibility.
\newblock {\em SIAM Journal on Matrix Analysis and Applications}, 19:487--498,
  1998.

\bibitem[Bar88]{Bar88}
N.~E Barabanov.
\newblock Lyapunov indicators of discrete linear inclusions, parts {I, II, and
  III}.
\newblock {\em Translation from Avtomat. e. Telemekh.}, 2, 3 and 5:40--46,
  24--29, 17--44, 1988.

\bibitem[Bar02]{Barvinok}
A.~Barvinok.
\newblock {\em A course in convexity}.
\newblock American Mathematical Society, 2002.

\bibitem[BN05]{BlNes05}
V.~D. Blondel and Yu. Nesterov.
\newblock Computationally efficient approximations of the joint spectral
  radius.
\newblock {\em SIAM J. Matrix Anal. Appl.}, 27(1):256--272, 2005.

\bibitem[BNT05]{BlNT04}
V.~D. Blondel, Yu. Nesterov, and J.~Theys.
\newblock On the accuracy of the ellipsoidal norm approximation of the joint
  spectral radius.
\newblock {\em Linear Algebra Appl.}, 394:91--107, 2005.

\bibitem[Bro74]{BrockettLie}
R.W. Brockett.
\newblock Lie algebras and {Lie} groups in control theory.
\newblock In D.Q. Mayne and R.W. Brockett, editors, {\em Geometric Methods in
  Systems Theory}, pages 17--56. D. Reidel Pub. Co., 1974.

\bibitem[BT80]{BraytonTong2}
R.~K. Brayton and C.~H. Tong.
\newblock Constructive stability and asymptotic stability of dynamical systems.
\newblock {\em IEEE Trans. Circuits and Systems}, 27(11):1121--1130, 1980.

\bibitem[BT00a]{BlTi2}
V.~D. Blondel and J.~N. Tsitsiklis.
\newblock The boundedness of all products of a pair of matrices is undecidable.
\newblock {\em Systems and Control Letters}, 41:135--140, 2000.

\bibitem[BT00b]{BlTi1}
V.~D. Blondel and J.~N. Tsitsiklis.
\newblock A survey of computational complexity results in systems and control.
\newblock {\em Automatica}, 36(9):1249--1274, 2000.

\bibitem[BW92]{BeWa92}
M.~Berger and Y.~Wang.
\newblock Bounded semigroups of matrices.
\newblock {\em Linear Algebra Appl.}, 166:21--27, 1992.

\bibitem[CH94]{CH94}
D.~Colella and C.~Heil.
\newblock Characterizations of scaling functions: continuous solutions.
\newblock {\em SIAM J. Matrix Anal. Appl.}, 15(2):496--518, 1994.

\bibitem[CLR95]{ChoiLamReznick}
M.-D. Choi, T.-Y. Lam, and B.~Reznick.
\newblock Sums of squares of real polynomials.
\newblock {\em Proceedings of Symposia in Pure Mathematics}, 58(2):103--126,
  1995.

\bibitem[DL92]{DaLa92}
I.~Daubechies and J.~C. Lagarias.
\newblock Sets of matrices all infinite products of which converge.
\newblock {\em Linear Algebra Appl.}, 161:227--263, 1992.

\bibitem[DL01]{DaLa01}
I.~Daubechies and J.~C. Lagarias.
\newblock Corrigendum/addendum to ``{Sets} of matrices all infinite products of
  which converge''.
\newblock {\em Linear Algebra Appl.}, 327:69--83, 2001.

\bibitem[DM99]{DayaMar}
W.~P. Dayawansa and C.~F. Martin.
\newblock A converse {Lyapunov} theorem for a class of dynamical systems which
  undergo switching.
\newblock {\em IEEE Transactions on Automatic Control}, 44:751--760, 1999.

\bibitem[GHND03]{GHNV}
Y.~Genin, Y.~Hachez, Yu. Nesterov, and P.~Van Dooren.
\newblock Optimization problems over positive pseudopolynomial matrices.
\newblock {\em SIAM J. Matrix Anal. Appl.}, 25(1):57--79 (electronic), 2003.

\bibitem[Gri96]{Gripenberg}
G.~Gripenberg.
\newblock Computing the joint spectral radius.
\newblock {\em Linear Algebra Appl.}, 234:43--60, 1996.

\bibitem[Joh48]{JohnEllipsoid}
F.~John.
\newblock Extremum problems with inequalities as subsidiary conditions.
\newblock In {\em Studies and Essays Presented to R. Courant on his 60th
  Birthday, January 8, 1948}, pages 187--204. Interscience Publishers, Inc.,
  New York, N. Y., 1948.

\bibitem[Koz90]{Koz90}
V.~A. Kozyakin.
\newblock Algebraic unsolvability of problem of absolute stability of
  desynchronized systems.
\newblock {\em Automation and Remote Control}, 51:754--759, 1990.

\bibitem[Lei92]{Leiz92}
A.~Leizarowitz.
\newblock On infinite products of stochastic matrices.
\newblock {\em Linear Algebra Appl.}, 168:189--219, 1992.

\bibitem[LP04]{LofbergParrilo}
J.~L\"ofberg and P.~A. Parrilo.
\newblock From coefficients to samples: a new approach to {SOS} optimization.
\newblock In {\em Proceedings of the 43$^{th}$ IEEE Conference on Decision and
  Control}, 2004.

\bibitem[Mae96]{Maesumi}
M.~Maesumi.
\newblock An efficient lower bound for the generalized spectral radius of a set
  of matrices.
\newblock {\em Linear Algebra Appl.}, 240:1--7, 1996.

\bibitem[Mar73]{MarcusMultilinear}
M.~Marcus.
\newblock {\em Finite dimensional multilinear algebra}.
\newblock M. Dekker, New York, 1973.

\bibitem[MM92]{MarcusMinc}
M.~Marcus and H.~Minc.
\newblock {\em A survey of matrix theory and matrix inequalities}.
\newblock Dover Publications Inc., New York, 1992.
\newblock Reprint of the 1969 edition.

\bibitem[Nes00]{NesterovSquared}
Yu. Nesterov.
\newblock Squared functional systems and optimization problems.
\newblock In {\em High performance optimization}, volume~33 of {\em Appl.
  Optim.}, pages 405--440. Kluwer Acad. Publ., Dordrecht, 2000.

\bibitem[NN94]{NN}
Y.~E. Nesterov and A.~Nemirovski.
\newblock {\em Interior point polynomial methods in convex programming},
  volume~13 of {\em Studies in Applied Mathematics}.
\newblock {SIAM}, {Philadelphia, PA}, 1994.

\bibitem[Par00]{Phd:Parrilo}
P.~A. Parrilo.
\newblock {\em Structured semidefinite programs and semialgebraic geometry
  methods in robustness and optimization}.
\newblock PhD thesis, California Institute of Technology, May 2000.
\newblock Available at
  \texttt{http://resolver.caltech.edu/CaltechETD:etd-05062004-055516}.

\bibitem[Par03]{sdprelax}
P.~A. Parrilo.
\newblock Semidefinite programming relaxations for semialgebraic problems.
\newblock {\em Math. Prog.}, 96(2, Ser. B):293--320, 2003.

\bibitem[PJ07]{PabloAliHSCC}
P.~A. Parrilo and A.~Jadbabaie.
\newblock Approximation of the joint spectral radius of a set of matrices using
  sum of squares.
\newblock In A.~Bemporad, A.~Bicchi, and G.~Buttazzo, editors, {\em Hybrid
  Systems: Computation and Control 2007}, volume 4416 of {\em Lecture Notes in
  Computer Science}, pages 444--458. Springer, 2007.

\bibitem[PK00]{ParriloKhatri}
P.~A. Parrilo and S.~Khatri.
\newblock On cone-invariant linear matrix inequalities.
\newblock {\em IEEE Transactions on Automatic Control}, 45(8):1558--1563, 2000.

\bibitem[Pro97]{Protasov1}
V.~Yu. Protasov.
\newblock The generalized joint spectral radius. {A} geometric approach.
\newblock {\em Izv. Ross. Akad. Nauk Ser. Mat.}, 61(5):99--136, 1997.
\newblock English translation in \emph{Izvestiya: Mathematics}, 61:5, 995-1030.

\bibitem[Pro05]{Protasov2}
V.~Yu. Protasov.
\newblock The geometric approach for computing the joint spectral radius.
\newblock In {\em Proceedings of the 44th IEEE Conference on Decision and
  Control and the European Control Conference 2005}, pages 3001--3006, 2005.

\bibitem[Rez00]{Reznick}
B.~Reznick.
\newblock Some concrete aspects of {H}ilbert's 17th problem.
\newblock In {\em Contemporary Mathematics}, volume 253, pages 251--272.
  American Mathematical Society, 2000.

\bibitem[RS60]{RoSt60}
G.~C. Rota and W.~G. Strang.
\newblock A note on the joint spectral radius.
\newblock {\em Indag. Math.}, 22:379--381, 1960.

\bibitem[RV06]{RohVandenberghe}
T.~Roh and L.~Vandenberghe.
\newblock Discrete transforms, semidefinite programming, and sum-of-squares
  representations of nonnegative polynomials.
\newblock {\em SIAM J. Optim.}, 16(4):939--964, 2006.

\bibitem[Sho87]{Shor}
N.~Z. Shor.
\newblock Class of global minimum bounds of polynomial functions.
\newblock {\em Cybernetics}, 23(6):731--734, 1987.
\newblock (Russian orig.: Kibernetika, No. 6, (1987), 9--11).

\bibitem[SWP97]{ShuWuPa97}
M.~Shih, J.~Wu, and C.~T. Pang.
\newblock Asymptotic stability and generalized {Gelfand} spectral radius
  formula.
\newblock {\em Linear Algebra Appl.}, 251:61--70, 1997.

\bibitem[TB97]{BlTi3}
J.~N. Tsitsiklis and V.D. Blondel.
\newblock The {L}yapunov exponent and joint spectral radius of pairs of
  matrices are hard- when not impossible- to compute and to approximate.
\newblock {\em Mathematics of Control, Signals, and Systems}, 10:31--40, 1997.

\bibitem[Tod01]{ToddSDP}
M.~Todd.
\newblock Semidefinite optimization.
\newblock {\em Acta Numerica}, 10:515--560, 2001.

\bibitem[VB96]{VaB:96}
L.~Vandenberghe and S.~Boyd.
\newblock Semidefinite programming.
\newblock {\em SIAM Review}, 38(1):49--95, March 1996.

\bibitem[Wir02]{wirth}
F.~Wirth.
\newblock Joint spectral radius and extremal norms.
\newblock {\em Linear Algebra Appl.}, 251:61--70, 2002.

\bibitem[WSV00]{HandSDP}
H.~Wolkowicz, R.~Saigal, and L.~Vandenberghe, editors.
\newblock {\em Handbook of Semidefinite Programming}.
\newblock Kluwer, 2000.

\bibitem[Zel94]{Zelen}
A.~L. Zelentsovsky.
\newblock Nonquadratic {L}yapunov functions for robust stability analysis of
  linear uncertain systems.
\newblock {\em IEEE Trans. Automat. Control}, 39(1):135--138, 1994.

\end{thebibliography}

\end{document}